\documentclass[a4paper,11pt]{scrartcl}
\usepackage[english]{babel}
\usepackage[utf8]{inputenc}
\usepackage[T1]{fontenc}
\usepackage{lmodern}
\usepackage[T1]{fontenc}
\usepackage{lmodern}
\usepackage[square,numbers]{natbib}
\usepackage{graphics}
\usepackage{graphicx}
\usepackage{amsmath}
\usepackage{latexsym}
\usepackage{siunitx}
\usepackage{amsfonts}
\usepackage{amssymb}
\usepackage[ngerman]{varioref}
\usepackage{hyperref}
\usepackage{cleveref}
\usepackage[abs]{overpic}
\usepackage{fullpage}
\usepackage[config=altsf]{subfig}
\begin{document}
\author{{\normalsize Karl Dieter Reinartz}}
\title{Chebychev \emph{inter}polations of the Gamma and Polygamma Functions and their analytical properties}
\maketitle
\textbf{in\,\,memoriam}\\
\textbf{Cornelius Lanczos\cite{Lanczos:1972} 1893-1974}\,\footnote{\,\,\,http://www.youtube.com/watch?v=avSHHi9QCjA\\
http://www.youtube.com/watch?v=PO6xtSxB5Vg}
\\\\
{\tiny address:}\\
{\tiny Email: KD.Reinartz@T-Online.de}\\
{\tiny Kieferndorfer Weg 30, D-91315 Höchstadt, GERMANY}\\
\\\textbf{keywords:}\\
Gamma function, Sterling formula, Bernoulli numbers, Chebychev approximations, Chebyshev polynomials, Shifted Chebyshev polynomials, Invers Gamma function,  Psi (Digamma) function, Harmonic function, Polygamma functions, Summation/Differentiation/Multiplication of Chebyshev approximations.\\
\tableofcontents
\section{Introduction}
The Gamma Function derived by Leonhard Euler (1729)\,is the generalization of discrete factorials:

\begin{equation}
z! = \Gamma  \left( z+1 \right) =\int _{0}^{\infty }\!{t}^{z-1}{e}^{-t}{dt},\quad \Re z>0
\end{equation}
The numerical evaluation is not easy. Whittaker+Watson \cite{Whittaker:1927:CMA} and Temme \cite{Temme:1996:SFI} give a good discussion of the $\Gamma$ Function and several basic properties.\\
In contrast to that Lanczos \cite{Lanczos:1964} developed approximations with a restricted precision by using Chebyshev polynomials in the range [-1..+1] (instead of shifted polynomials which will be used exclusively in this paper).
Chebychev polynomials were introduced into numerical analysis especially by Lanczos \cite{Lanczos:1972}in the US since 1935 and by Clenshaw \cite{Clenshaw:1962:CSM} in GB since 1960.--\\\\
The next formula is due to James Stirling (1730)
\begin{equation}
\label{eqn:lngamma sterling}
\ln   \Gamma  \left( z \right)  \,\sim\,ln  \left( \sqrt {2\pi}\,
{z}^{z-\frac{1}{2}}\,{{\rm e}^{-z}} \right) +\sum _{n=1}^{\infty }{\frac {B_{{2\,n}}}{2\,n \left( 2\,n-1 \right) {z}^{2\,n-1}}}
\end{equation}
The $B_{{2\,n}}$ are the Bernoulli numbers with a poor behaviour:
\begin{equation}
\begin{split}
&\frac{1}{6},\,-\frac{1}{30},\,\frac{1}{42},\,-\frac{1}{30},\,\frac {5}{66},\,\frac {691}{2730},\,\frac{7}{6},\,-\frac{3617}{510},\,\frac{43867}{798},\, -\frac{174611}{330},\,\frac {854513}{138},-\frac {236364091}{2730},\\&\quad\frac{8553103}{6},\,-\frac{23749461029}{870},\,\frac{8615841276005}{14322},\,\frac{7709321041217}{510},\,\frac {2577687858367}{6},...
\end{split}
\end{equation}
They decrease at the beginning only slowly and then grow with (2n)!.
The complete terms in the infinite sum 
\cref{eqn:lngamma sterling} 
depend on n and z and grow nevertheless, especially if z is small:
\begin{equation}
\begin{split}
&\frac{1}{12 \cdot z},\,-\frac {1}{360 \cdot z^3},\,
\frac{1}{1260 \cdot z^5},\,-\frac{1}{1680 \cdot z^7},\,\frac{1}{1188 \cdot{
z^9}},\,-\frac {691}{360360 \cdot z^{11}},\,\frac {1}{156 \cdot z^{13 }},\,
-\frac {3617}{122400 \cdot z^{15}},\\
&\quad \frac {43867}{244188 \cdot z^{17 }},-\frac {174611}{125400 \cdot z^{19 }},\,
\frac {77683}{5796 \cdot z^{21 }},\,-\frac {236364091}{1506960 \cdot z^{23 }},\,\frac {657931}{300 \cdot z^{25 }},\,-\frac {3392780147}{93960 \cdot z^{27 }},\\
&\quad \frac {1723168255201}{2492028 \cdot z^{29 }},\,
-\frac {7709321041217}{505920 \cdot z^{31 }},\frac {151628697551}{396 \cdot z^{33}},\,-\frac {26315271553053477373}{
2418179400 \cdot z^{35 }},...\\
\end{split}
\end{equation}
\begin{figure}[!ht]
\centering \unitlength .5cm
\subfigure[summation limit]
{\begin{picture}(12,13)
\includegraphics[angle=0,scale=0.3]{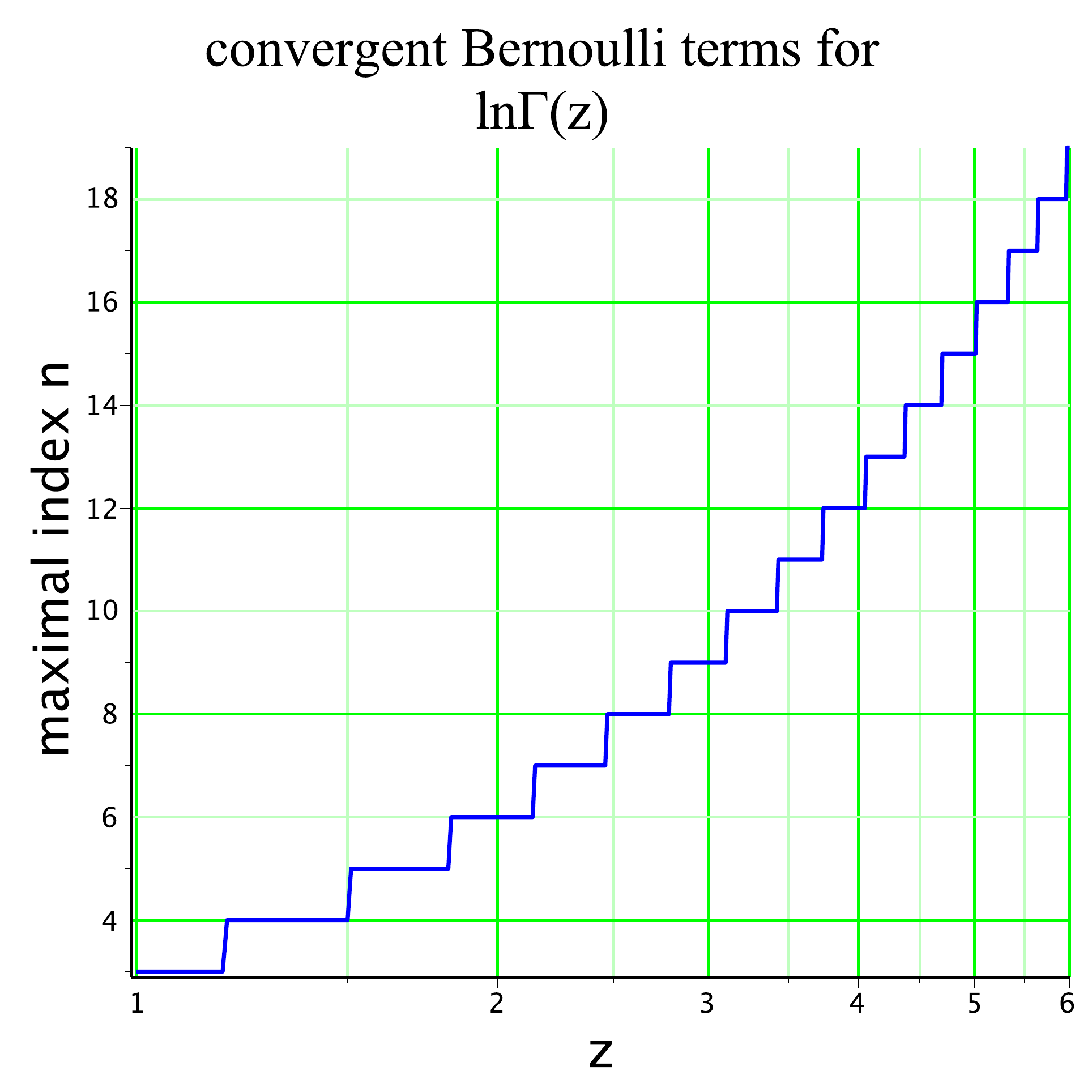}
\label{fig:lnGconv}
\end{picture}}
\qquad
\subfigure[correct decimal digits]
{\begin{picture}(12,13)
\includegraphics[angle=0,scale=0.3]{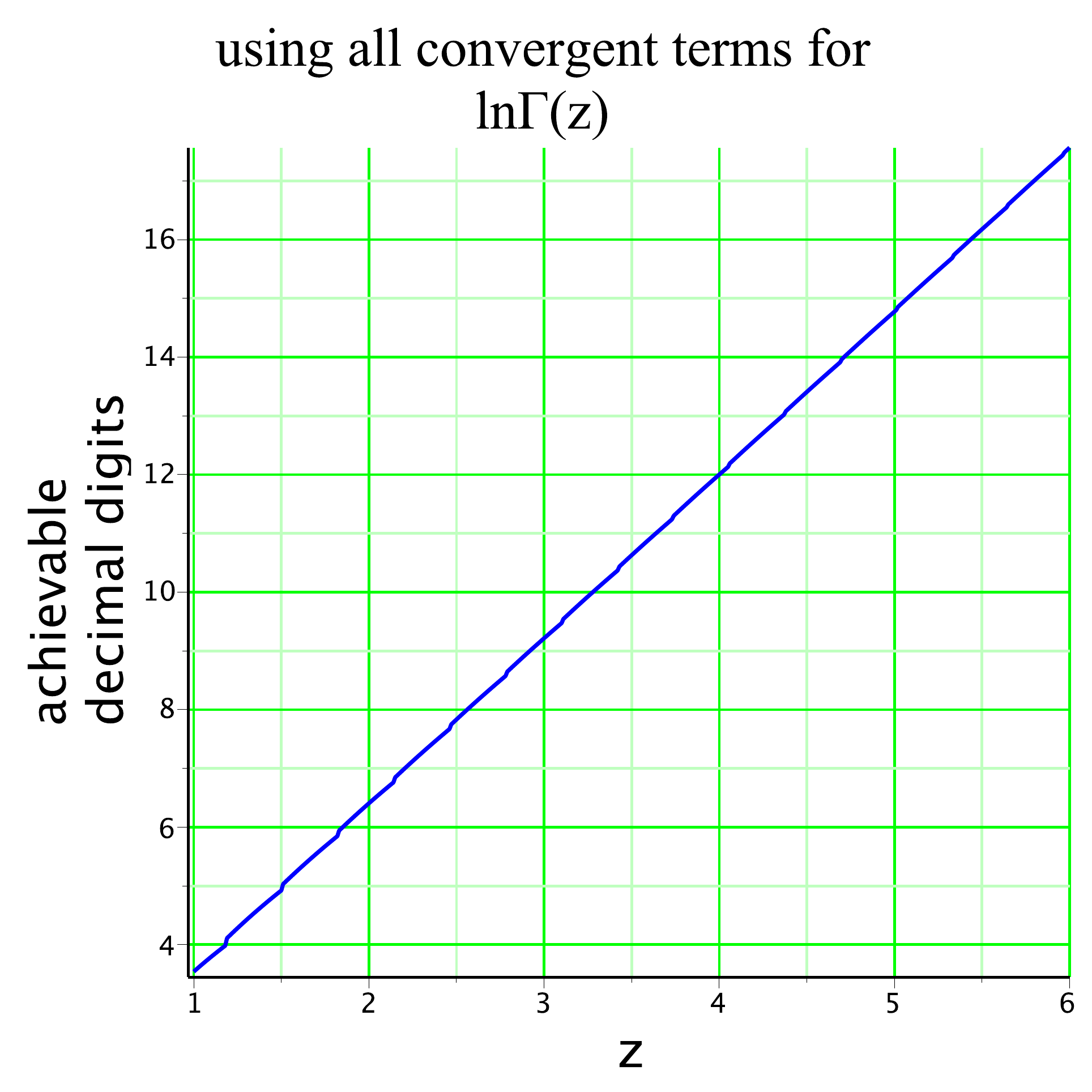}
\label{fig:lnGbest}
\end{picture}}
\caption{restricted summation and limited precision}
\label{fig:total}
\end{figure}

The summation has to stop \textbf{before the terms begin to grow unrestricted}. There is an optimal position depending on z where summation has to end. This problem is discussed in some detail in\cite{Graham:1994:CM}\footnote{...page 467} .\\
A further disadvantage is the low convergence of the admitted terms. In fig. \ref{fig:total}
the problem is described in some detail depending on z:
fig. \ref{fig:lnGconv} shows the maximal number of convergent terms,
fig. \ref{fig:lnGbest} shows the maximal achievable accuracy in decimal digits.\\
\section{Chebyshev \emph{Inter}polations of the $\Gamma$\,Function}
The Chebyshev polynomials were derived by the Russian Mathematician P. L. Chebyshev (1821-1894) \cite{Natanson:1955}. Among all normalized power polynomials of same degree they have the smallest deviation from zero in a predefined intervall. Most of their beautiful properties are described by Snyder \cite{Snyder:1966} and Clenshaw \cite{Clenshaw:1962:CSM} showing many applications to transcendental functions and differential equations.\\
The approximation of the $\Gamma$ Function is represented by
\begin{equation}
(z-1)!\, \simeq  \Gamma  \left( z \right)  = \sqrt {2\pi}\,
{z}^{z-\frac{1}{2}}\,{{\rm e}^{-z}} \, *\,\sideset{_{\,0}^{\infty}}{_{}^{'}}\sum a{_r^*}T{_r^*}(\frac{1}{z})\, ,\quad 1\leq ~z~ \leq \infty
\end{equation}
\begin{figure}[!ht]
\centering
\includegraphics[angle=90,scale=1.25]{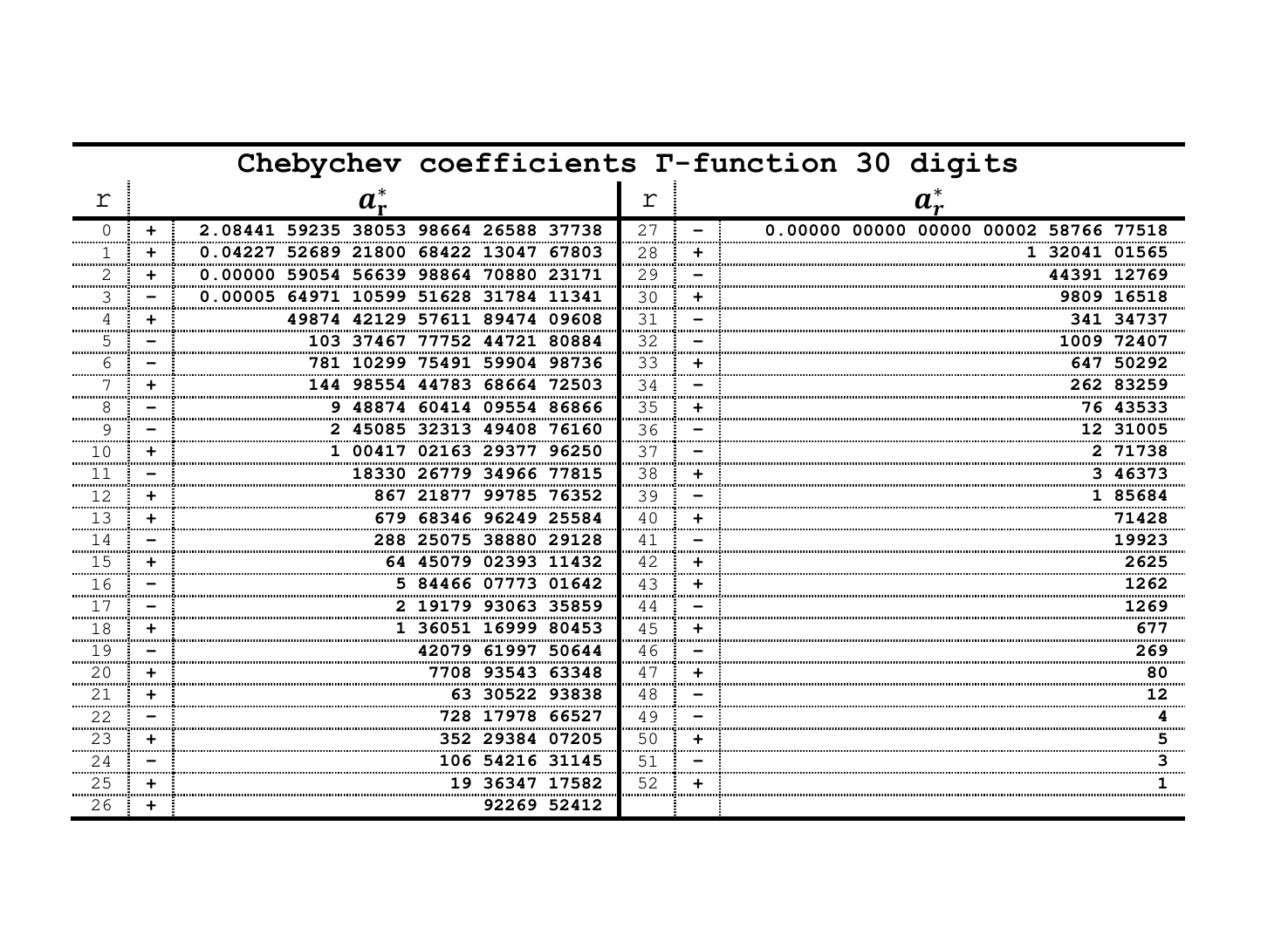}
\caption{Chebyshev coefficients $\Gamma$-function 30 digits}
\label{fig:GAMMA}
\end{figure}
\Cref{fig:GAMMA}\,\,contains 53 Chebyshev coefficients of the $\Gamma$-function for an accuracy of 30 decimal digits in the whole range $1\leq z \leq \infty$.\,Using two coefficients the corresponding powerseries is:
\begin{figure}[!ht]
\centering \unitlength .5cm
\subfigure[relative error using only 2 coefficients]
{\begin{picture}(13,12)
\includegraphics[angle=0,scale=0.3]{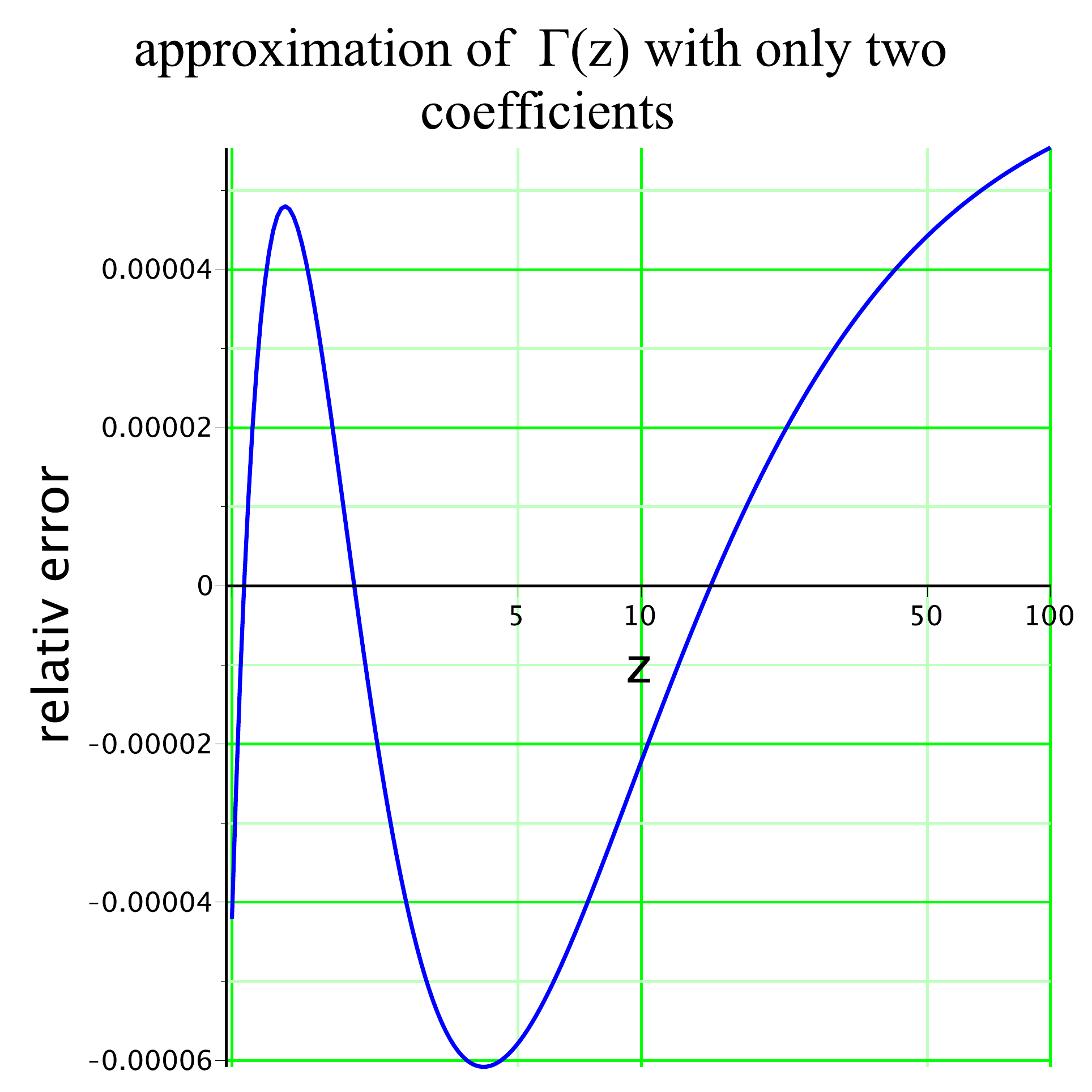}
\label{fig:Gconv1}
\end{picture}}
\qquad
\subfigure[relative error using 11 coefficients]
{\begin{picture}(13,12)
\includegraphics[angle=0,scale=0.3]{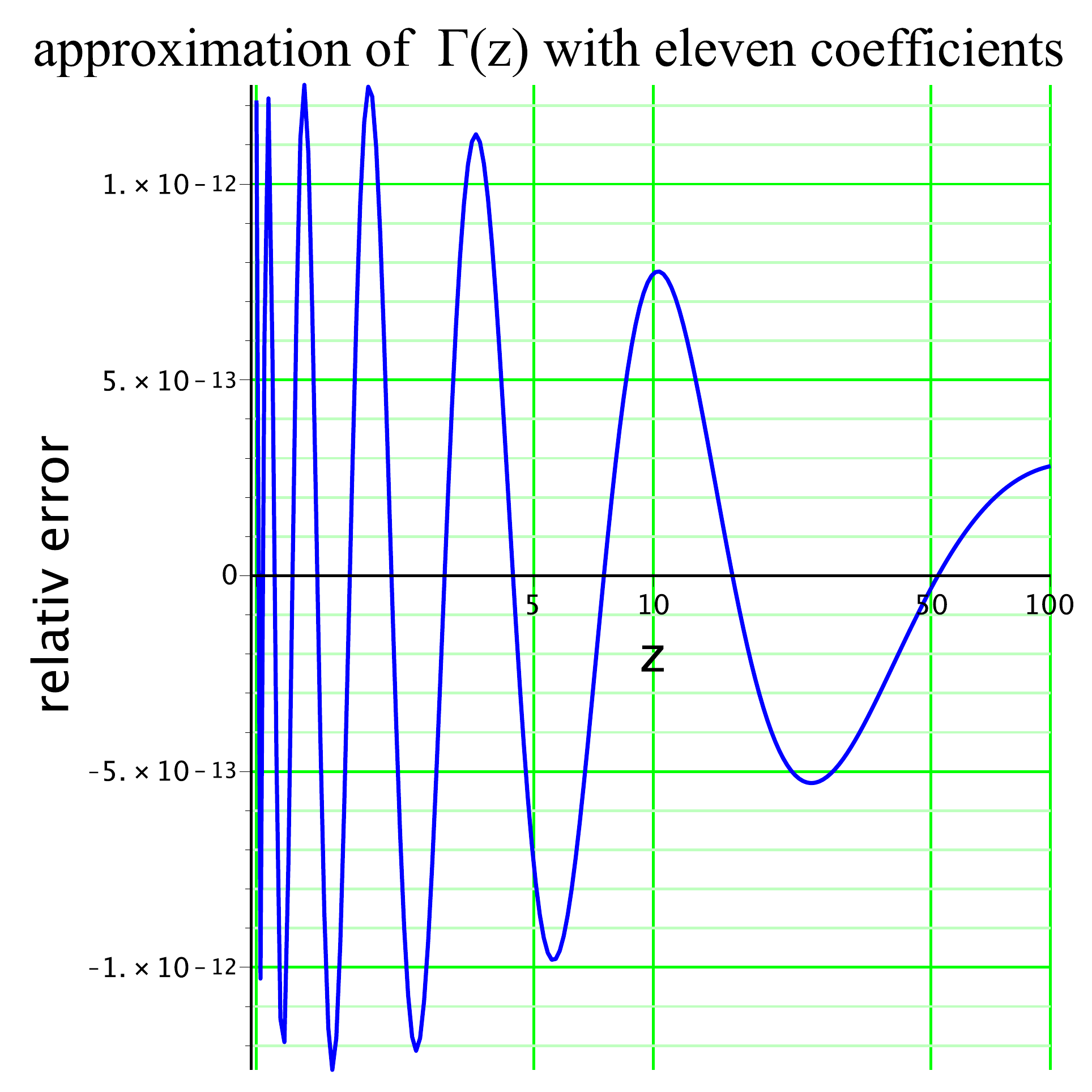}
\label{fig:Gconv2}
\end{picture}}
\caption{overall convergence}
\label{fig:Gtotal}
\end{figure}
\begin{equation}
\begin{split}
\Gamma  \left( z \right) & = \sqrt {2\pi}\cdot
z^{z-\frac{1}{2}}\cdot{\rm e}^{-z} * [0.999935+\frac{0.0845506}{z}]\\
\end{split}
\end{equation}
the maximal relativ error(\Cref{fig:Gconv1}) is less than $8*10^{-4}$. \.
Using eleven coefficients for the corresponding powerseries
\begin{equation}
\begin{split}
\Gamma  \left( z \right) & = \sqrt {2\pi}\cdot
z^{z-\frac{1}{2}}\cdot{\rm e}^{-z} \\ 
& \quad * [ 0.99999999998+\frac{0.083333337647}{z}+\frac{0.0034720552506}{z^2}-\frac{0.0026788696285}{z^3}\\
& \quad +\frac{0.00024711193390}{z^4}+\frac{0.00084986066787}{z^5}-\frac{0.000035855790507}{z^6}-\frac{0.00068599470338}{z^7}\\
& \quad +\frac{0.00067284352663}{z^8}-\frac{0.00029536102066}{z^9}+
\frac{0.000052647439438}{z^{10}} ]
\end{split}
\end{equation}
the maximal relativ error(\Cref{fig:Gconv2}) is less than $2*10^{-11}$.
\subsection{The $\Gamma^{-1}$ Function}
\begin{figure}[htbp]
\centering
\includegraphics[angle=90,scale=1.5]{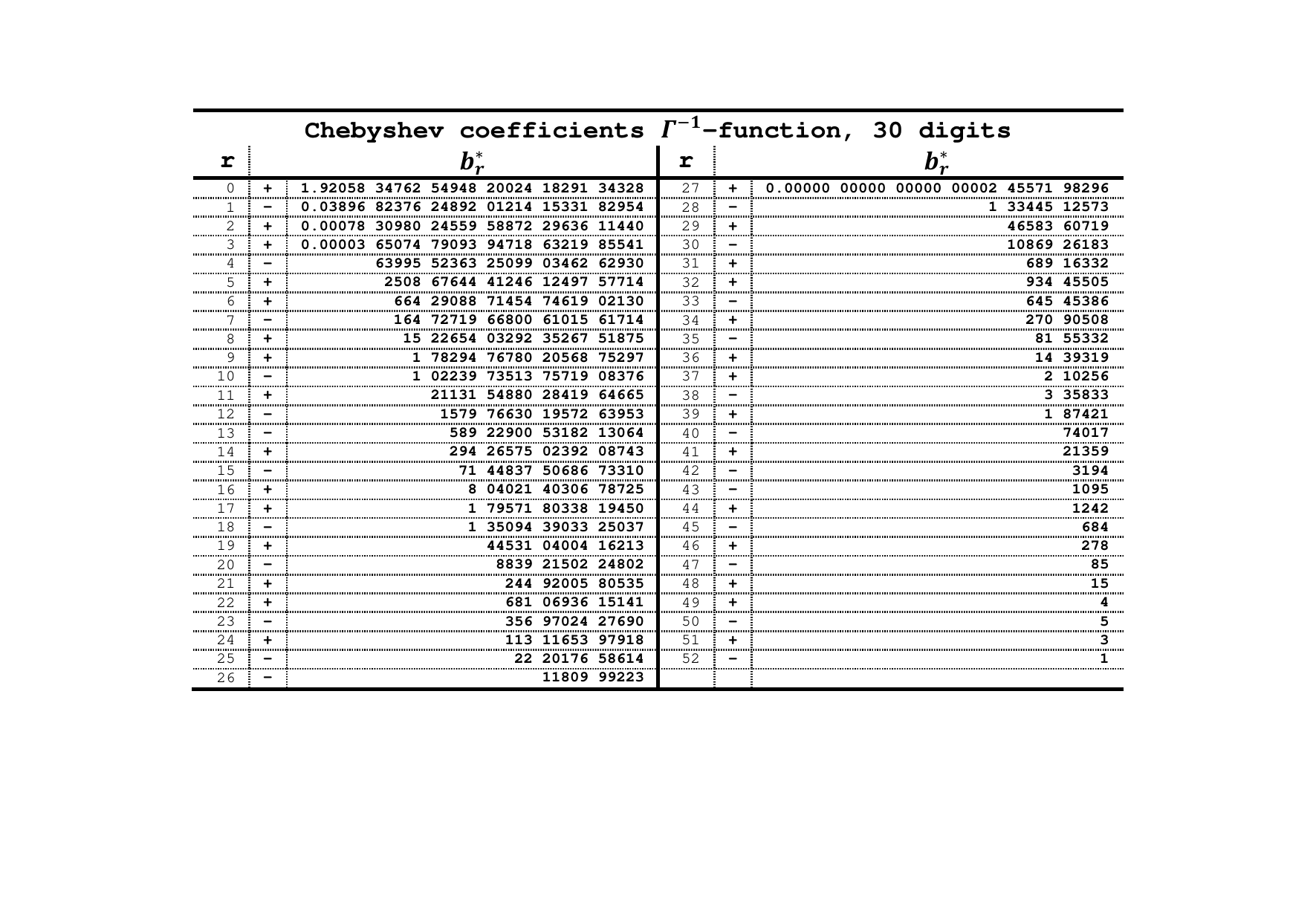}
\caption{Chebyshev coefficients $\Gamma^{-1}$-function 30 digits}
\label{fig:invGAMMA}
\end{figure}
\begin{figure}[htbp]
\centering
\includegraphics[scale=0.4]{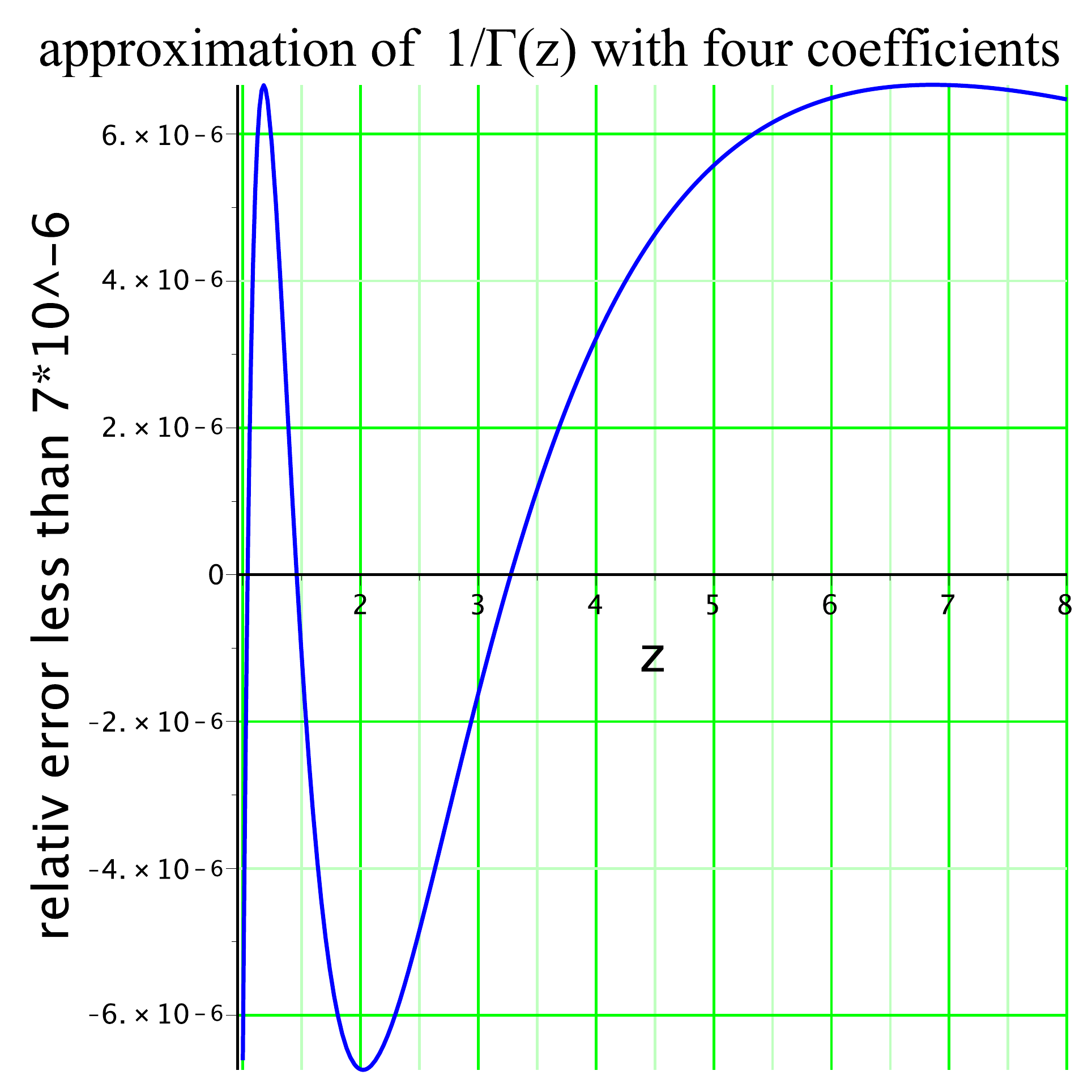}
\caption{ $\Gamma^{-1}$-function error analysis}
\label{fig:invGAMMAerr}
\end{figure}
\Cref{fig:invGAMMA}\,contains 53 Chebyshev coefficients of the $\Gamma^{-1}$-function for an accuracy of 30 decimal digits in the whole range $1\leq z \leq \infty$.
\begin{equation}
\frac{1}{(z-1)!}\, \simeq  \Gamma^{-1}  \left( z \right)  
=\frac{1}{\sqrt {2\pi}}\,
{z}^{\frac{1}{2}-z}\,{{\rm e}^{z}} \, *\,\sideset{_{\,0}^{\infty}}{_{}^{'}}\sum b{_r^*}T{_r^*}(\frac{1}{z})\, ,\quad 1\leq ~z~ \leq \infty
\end{equation}\\
With four coefficients the powerseries expansion is:\\
\begin{equation}
\Gamma^{-1}  \left( z \right)  
=\frac{1}{\sqrt {2\pi}}\,
{z}^{\frac{1}{2}-z}\,{{\rm e}^{z}} \, *\, [  1.000006-
 \frac{0.08354413}{z}\,+\frac{0.004512425}{z^2} \,+\frac{0.001168239}{z^3}] 
\end{equation}\\
The maximal relativ error (\Cref{fig:invGAMMAerr}) is less than $7*10^{-6}$.\\
In contrast to that in the famous Handbook of Mathematical Functions \cite{Abramowitz:1964:HMF} \footnote{...page 256} the series expansion for $\Gamma^{-1}$ is completely wrong.
\subsection{The Ln$\Gamma$ Function}
\begin{figure}[!ht]
\centering
\includegraphics[angle=90,scale=1.25]{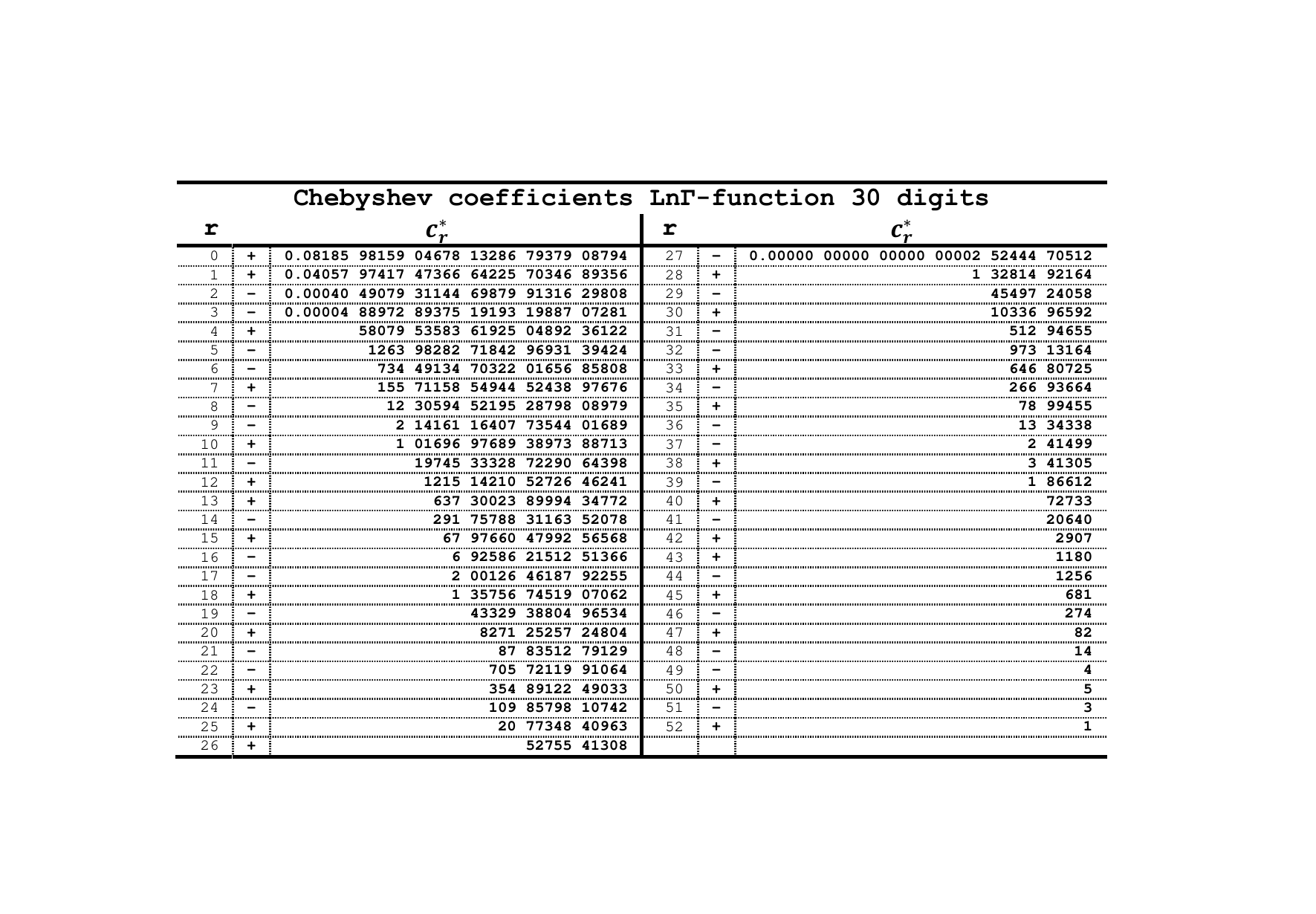}
\caption{Chebyshev coefficients $ln\Gamma$-function 30 digits}
\label{fig:lnGAMMA}
\end{figure}
\begin{figure}[!ht]
\centering \unitlength .5cm
\subfigure[absolute error using only 2 coefficients]
{\begin{picture}(13,12)
\includegraphics[angle=0,scale=0.3]{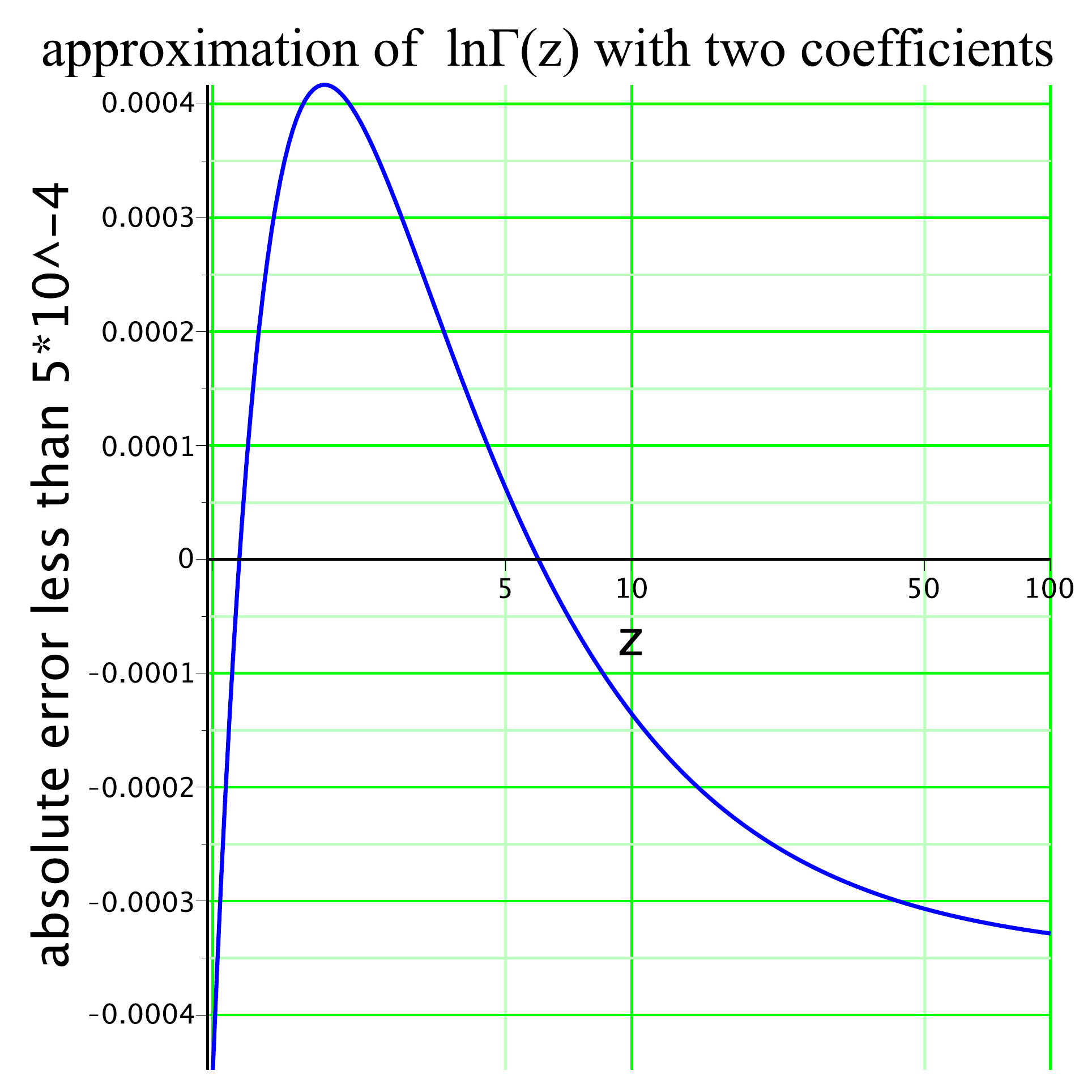}
\label{fig:lnGconv1}
\end{picture}}
\qquad
\subfigure[absolute error using 5 coefficients]
{\begin{picture}(13,12)
\includegraphics[angle=0,scale=0.3]{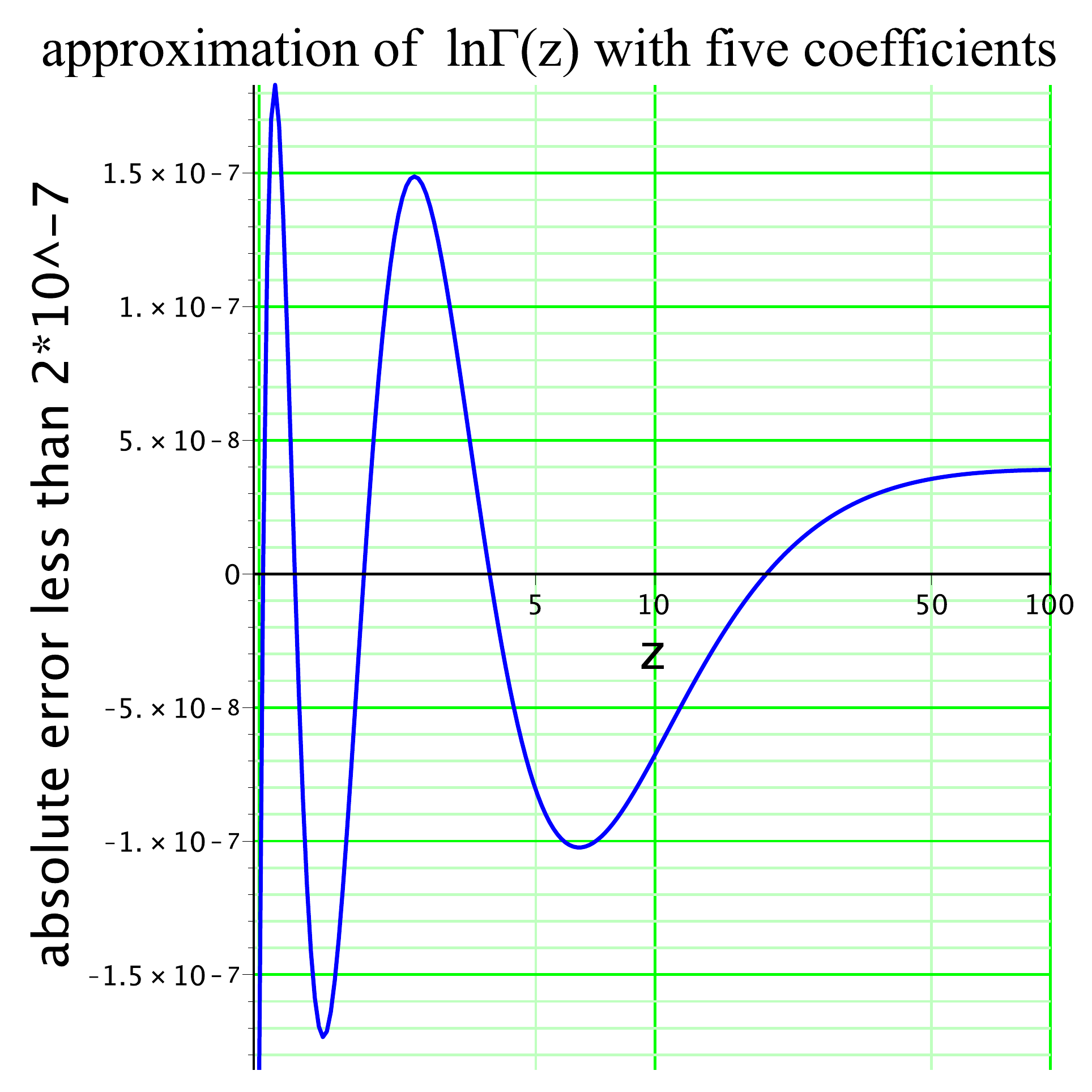}
\label{fig:lnGconv2}
\end{picture}}
\caption{overall convergence}
\label{fig:lnGtotal}
\end{figure}
The approximation is represented by
\begin{equation}
\ln[(z-1)!]\, \simeq \ln   \Gamma  \left( z \right)  =\ln  \left( \sqrt {2\pi}\,
{z}^{z-\frac{1}{2}}\,{{\rm e}^{-z}} \right) +\sideset{_{\,0}^{\infty}}{_{}^{'}}\sum c{_r^*}T{_r^*}(\frac{1}{z})\, ,\quad 1\leq ~z~ \leq \infty
\end{equation}
\Cref{fig:lnGAMMA} contains the Chebyshev coefficients with a precision of 30 decimal digits for the whole range of $1\leq ~z~ \leq \infty$. Using only the first two coefficients and building the power series form
\begin{equation}
\ln   \Gamma  \left( z \right)  =\ln   \sqrt {2\pi} \,+  \left( z- \frac{1}{2} \right) \ln  \left( z \right)- z + 0.91932+ \frac{0.081160}{z}
\end{equation}
the maximum absolute error (\Cref{fig:lnGconv1}) is less than $5*10^{-4}$.
Using five coefficients
\begin{equation}
\begin{split}
\ln   \Gamma  &\left( z \right)  = \ln   \sqrt {2\pi}\,+ \left( z- \frac{1}{2} \right) \ln   z \,- z \\ & \quad
+ 0.918935\, + \frac{0.0833326}{z}\,+\frac{0.000037082}{z^{2}} -\frac{0.00305155}{z^{3}}\,+ 
\frac{0.000743418}{z^{4}}\\
\end{split}
\end{equation}
the maximum absolute error (\Cref{fig:lnGconv2}) is less than $2*10^{-7}$.\\
\section{The Chebyshev \emph{Inter}polation of the Psi\,(Digamma) Function}
This function is the first derivative of the Ln$\Gamma$ Function:
\begin{equation}
\psi^{(0)}(z)\,=\,\psi(z)=\ln^{'}   \Gamma  \left( z \right)  =\frac{\Gamma ^{'}(z)}{\Gamma(z)}=ln\,z-\frac{1}{2z} +\sideset{_{\,0}^{\infty}}{_{}^{'}}\sum c{_r^*}T{_r^{*'}}(\frac{1}{z})\, ,\quad 1\leq ~z~ \leq \infty
\end{equation}
After differentiating the sum using \cref{eqn:diff} and \cref{eqn:diff+} , $-\frac{1}{2z}=-\frac{1}{4}*
(T^{*}_{0}(\frac{1}{z})+T^{*}_{1}{(\frac{1}{z})})$ \,has to be added. The final result is
\begin{equation}
\label{eqn:Psi}
\psi^{(0)}(z)\,=\,\psi(z) =\frac{\Gamma ^{'}(z)}{\Gamma(z)}=ln\,z +\sideset{_{\,0}^{\infty}}{_{}^{'}}\sum \alpha{_{(0)r}^*}T{_r^{*}}(\frac{1}{z})\, ,\quad 1\leq ~z~ \leq \infty
\end{equation}
\begin{figure}[htbp]
\centering
\includegraphics[scale=1.3]{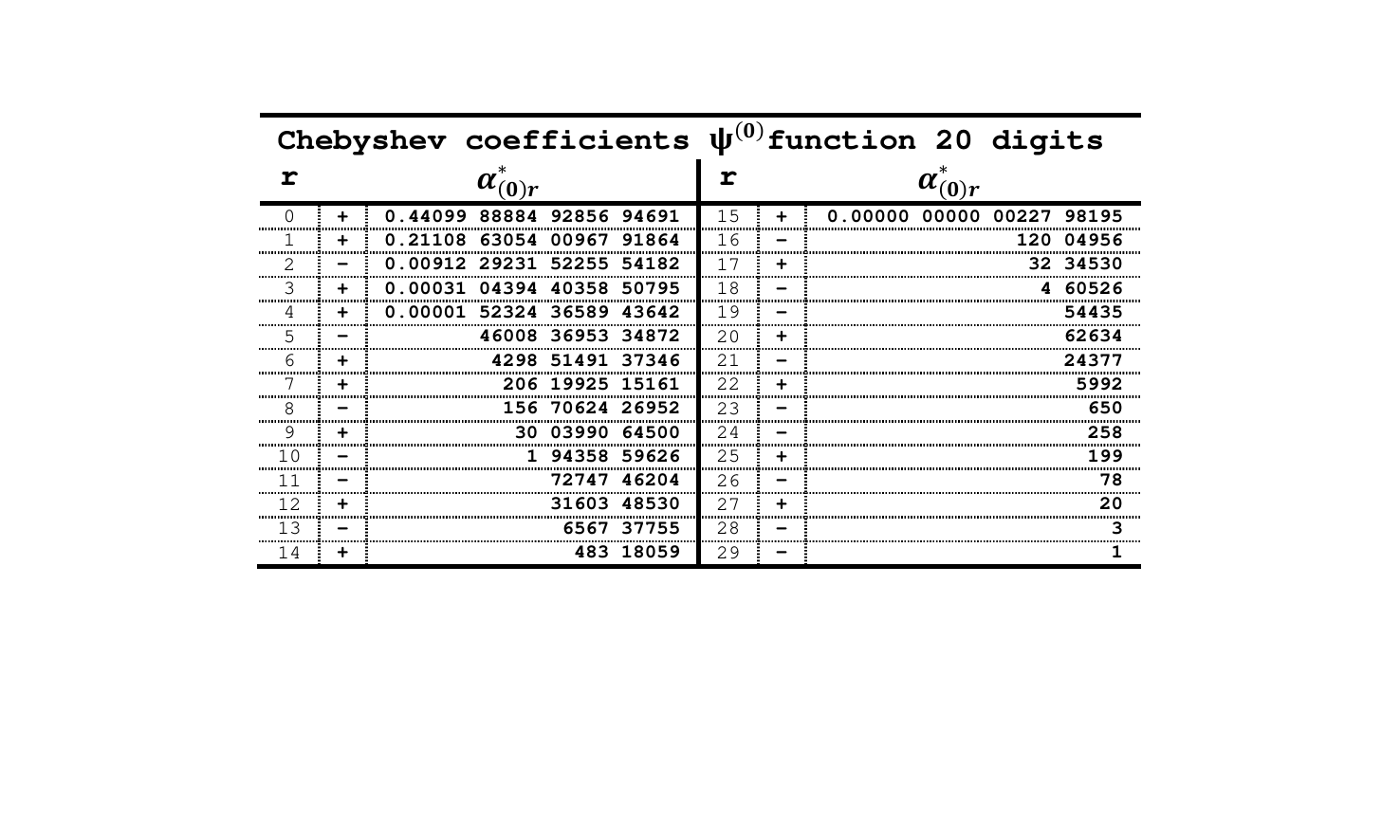}
\caption{The $\psi^{(0)}$-approximation 20 digits}
\label{fig:psio}
\end{figure}
\\
\subsection{Summation of the Harmonic Series}
$-\psi^{(0)}(1)\,=\,\gamma\,= 0.57721\,\,56649\,\,01532\,\,86061\,\,$is Euler's constant. \\\\
$\psi^{(0)}(n+1)\,-\,\psi^{(0)}(1)\,=\,H_{n}\,=\,1+\frac{1}{2}+\frac{1}{3}+\,...+\frac{1}{n}\,,\quad n\in N$\\\\defines and computes the n-th harmonic number $H_{n}$.
\section{\emph{Inter}polating further Polygamma Functions}
Differentiating the result of \cref{eqn:Psi} as before one gets
\begin{equation}
\psi^{(1)}(z) =\psi^{'}(z) =\frac{1}{z} +\sideset{_{\,0}^{\infty}}{_{}^{'}}\sum \alpha{_{(0)r}^*}T{_r^{*'}}(\frac{1}{z})\, ,\quad 1\leq ~z~ \leq \infty
\end{equation}
Finally $\frac{1}{z}=\frac{1}{2}*
(T^{*}_{0}({\frac{1}{z}})+T^{*}_{1}({\frac{1}{z}}))$\,has to be added yielding
\begin{equation}
\psi^{(1)}(z) =\sideset{_{\,0}^{\infty}}{_{}^{'}}\sum \alpha_{(1)r}^{*}T{_r^{*}}(\frac{1}{z})\, ,\quad 1\leq ~z~ \leq \infty
\end{equation}
The higher Polygamma Functions can be approximated applying the two step differentiation repeatedly without additional correction. Each next generated function looses about two decimal digits in precision.
\subsection{Summation of the higher Harmonic Series}
The general relation is:
\begin{equation}
\frac{(-1)^{m+1}}{m!}\,\psi^{(m)}(z) =\sum{_{k=\,0}^{\infty}} \frac{1}{(z+k)^{m+1}}\,
=\,\frac{1}{z^{m+1}}+\frac{1}{(z+1)^{m+1}}+\frac{1}{(z+2)^{m+1}}+...
\end{equation}
and especially for z=n integer
\begin{equation}
\frac{(-1)^{m+1}}{m!}\, [\psi^{(m)}(1)-\psi^{(m)}(n)]\,=\,\frac{1}{1^{m+1}}+\frac{1}{2^{m+1}}+\frac{1}{3^{m+1}}+\,...\,+\frac{1}{(n-1)^{m+1}}
\end{equation}
and further specialized with m=1
\begin{equation}
\psi^{(1)}(1)-\psi^{(1)}(n)\,=\,\frac{1}{1^{2}}+\frac{1}{2^{2}}+\frac{1}{3^{2}}+\,...\,+\frac{1}{(n-1)^{2}}
\end{equation}
\begin{figure}[htbp]
\centering
\includegraphics[scale=1.2]{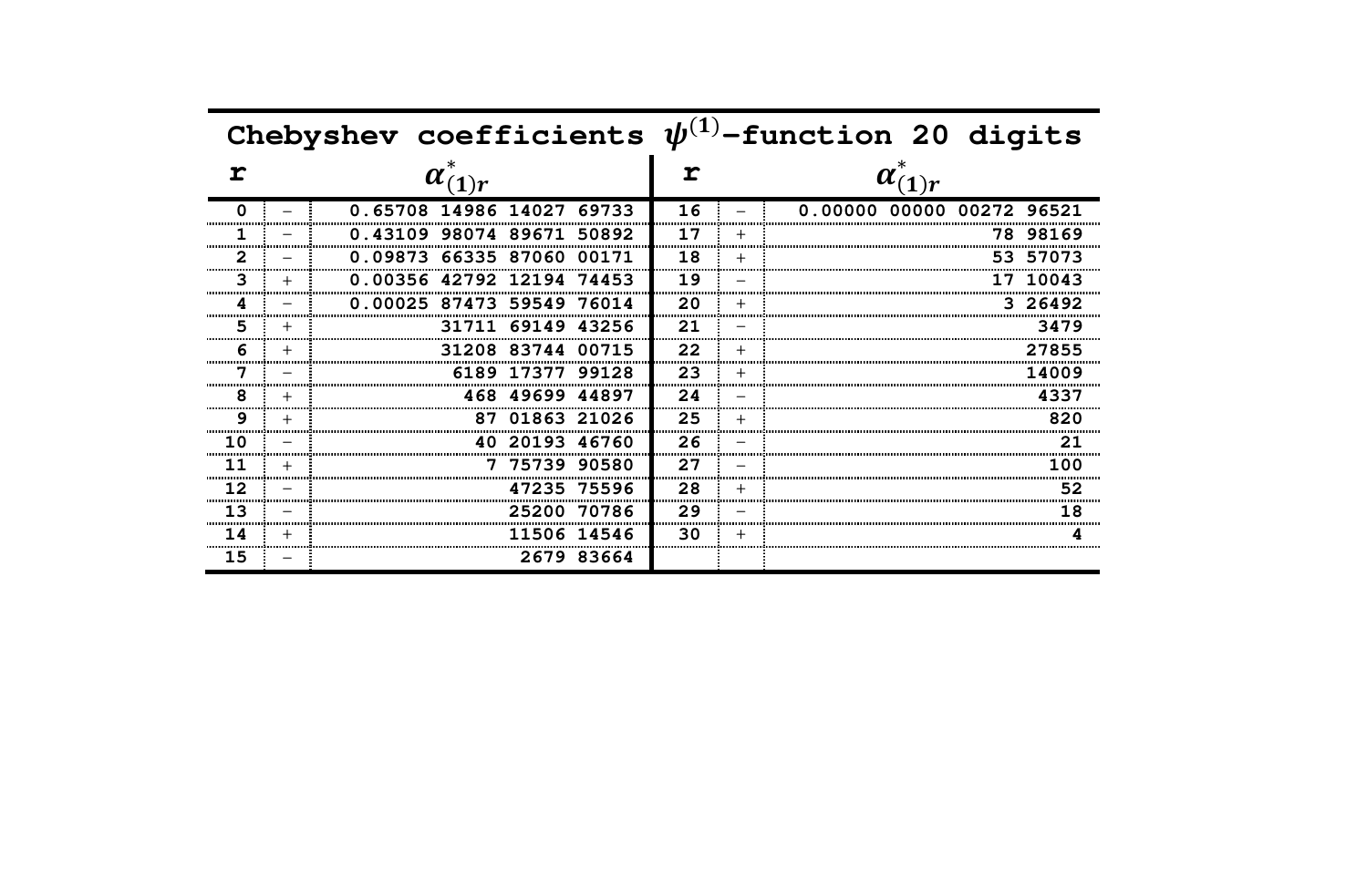}
\caption{The $\psi^{(1)}$-function 20 digits}
\label{fig:psi1}
\end{figure}
\section{Relations of the Shifted Chebychev Polynomials}
The Shifted Chebyshev polynomials are defined by
\begin{align*}
T^{*}_{r}(\theta)&=cos(r\theta) \qquad -1\leq\cos\,\theta\leq+1 \quad &-1\leq T^{*}_{r}\leq+1 \\
2z-1&=cos\,\theta\quad &0\leq ~z~ \leq+1
\end{align*}
They are power polynomials in z. Their highest coefficient $2^{r-1}$ is used for normalization.\\
Chebyshev proved \cite{Natanson:1955} that among all normalized power polynomials of same degree (or less) they have the smallest deviation from zero in the range $0\leq z\leq+1$. That makes them unique for optimal interpolation in the declared region. The ranges may be adapted by linear or even nonlinear transformations.\\
The polynomials for the intervall $0\leq z\leq +1$ are called the shifted polynomials. They are used here exclusively.
Explicit expressions for the first few shifted Chebyshev polynomials are:
$T^{*}_0(z)=1,\quad T^{*}_1(z)=2z-1,\quad T^{*}_2(z)=8z^{2}-8z+1,\quad T^{*}_3(z)=32z^3-48z^2+18z-1,...$\\\\Inversion gives:\\
$1=T^{*}_0(z),\quad 2z=T^{*}_0(z)+T^{*}_1(z),\quad 8z^2=3T^{*}_0(z)+4T^{*}_1(z)+T^{*}_2(z), \quad \\32z^{3}=10\,T^{*}_{0}(z)+15\,T^{*}_{1}(z)+6\,T^{*}_{2}(z)+T^{*}_{3}(z),$...\\
\subsection{Chebychev approximation of smooth functions}
\begin{equation}
\label{eqn:chebfunct}
f(z)=\sideset{_0^{n}}{_{}^{'}}\sum a{_r^*}T{_r^*}(z)=\frac{1}{2}a{^*_0}
+a{^*_1}T{^*_1}(z)+a{^*_2}T{^*_2}(z)+a{^*_3}T{^*_3}(z)+...
\end{equation}
\subsubsection{Numerical determination of the $a{_r^*}$-coefficients}
For the given function f(z) the $a{_r^*}$ can be determined by
\begin{equation}
\label{eqn:coeffsr}
a{_r^*}=\sideset{_{j=0}^{m}}{_{}^{''}}\sum f(cos^2(\frac{j\pi}{2m}))\,
cos(\frac{rj\pi}{m})
\end{equation}
$^{''}$\,means: terms with j=0 and j=m must be halfed and
m should be chosen sufficiently large\\
for a good approximation.
\subsubsection{Summation}
\begin{enumerate}
\item substituting the $T^{*}$-polynomials by their powerseries representations and thereafter applying the Horner Scheme or
\item it is better to use the coefficients directly: starting with a sufficiently large index n and
applying recursion:\\
\label{eqn:sum}
\begin{equation}
b^{*}_{r} = (2*z-1)*b^{*}_{r+1}-b^{*}_{r+2}+a^{*}_{r}\,,\quad b^{*}_{n+1} = b^{*}_{n+2}=0\,,\quad  r=n,n-1,...,0
\end{equation}
\center $f(z)=\frac{1}{2}(b^{*}_{0}-b^{*}_{2})$\\
\end{enumerate}
\subsubsection{Differentiation}
\begin{enumerate}
\item In order to get $f^{'}(z)=\sideset{_0^{n-1}}{_{}^{'}}\sum a{_r^{*'}}T{_r^*}(z)$ from \cref{eqn:chebfunct} one starts with a sufficiently large index r=n
\begin{equation}
\label{eqn:diff}
a^{*'}_{r-1}= a^{*'}_{r+1}+4\,r\,a^{*}_{r}, \quad a^{*'}_{n}=0, \quad a^{*'}_{n+1}=0
\end{equation}
and applies the recursion till r=1.
\item Differentiation (chainrule) of \center $f(x)=\sideset{_{\,0}^{n}}{_{}^{'}}\sum a{_r^*}T{_r^*}({x})$with $x=\frac{1}{z}$ results in
$$f^{'}(z)=-\frac{1}{z^{2}}*\sideset{^{n-1}_{0}}{_{}^{'}}\sum a{_r^{*'}}T{_r^*}(\frac{1}{z})$$\\
In addition to the former derivation step each coefficient of the derived form has to be multiplied by \\ 
\begin{equation}
\label{eqn:diff+}
-\frac{1}{z^{2}}=-(\frac{3}{8}T^{*}_{0}(\frac{1}{z})+\frac{1}{2}T^{*}_{1}(\frac{1}{z})
+\frac{1}{8}T^{*}_{2}(\frac{1}{z}))
\end{equation}
\center applying the multiplication rule of the next subsection.
\end{enumerate}
\subsubsection{Multiplication of two Chebyshev approximations}
The relation\\
\begin{equation}
T^{*}_{m}{(z)}*T^{*}_{n}{(z)} = \frac{1}{2}*[T^{*}_{m+n}{(z)}+ T^{*}_{|m-n|}{(z)}]
\end{equation}\\ 

is used for multiplying two polynomials \cref{eqn:chebfunct}. The resulting polynomial has m+n+1 coefficients and may be further reduced in length with a minor loss in accuracy.\\\\
\bibliography{literatur_21}
\bibliographystyle{natdin}
\end{document}